\documentclass[hidelinks,onefignum,onetabnum]{siamart251216}
\usepackage{lipsum}
\usepackage{amsfonts}
\usepackage{graphicx}
\usepackage{epstopdf}
\ifpdf
  \DeclareGraphicsExtensions{.eps,.pdf,.png,.jpg}
\else
  \DeclareGraphicsExtensions{.eps}
\fi

\newsiamremark{remark}{Remark}
\newsiamremark{hypothesis}{Hypothesis}
\crefname{hypothesis}{Hypothesis}{Hypotheses}
\newsiamthm{claim}{Claim}
\newsiamremark{fact}{Fact}
\crefname{fact}{Fact}{Facts}

\headers{Trust-region high-index saddle dynamics}{G. Ji, L. Li, Z. Xu, and J. Yin}

\title{
Trust-region high-index saddle dynamics for ill-conditioned saddle-point computation
\thanks{This work was partially supported by the
National Natural Science Foundation of China
(Grants No.~12471363, No.~12201021 and No.~12601819).%
}}

\author{
Guanghua Ji%
\thanks{School of Mathematical Sciences,
Laboratory of Mathematics and Complex Systems,
Ministry of Education, 
Beijing Normal University,
Beijing 100875, China
(\email{ghji@bnu.edu.cn},
\email{202321130105@mail.bnu.edu.cn},
\email{jyyin@bnu.edu.cn}).}
\and
Longjing Li\footnotemark[2]
\thanks{Corresponding author
(\email{202321130105@mail.bnu.edu.cn}).}
\and
Zhen Xu%
\thanks{School of Mathematics, Statistics and Mechanics,
Beijing University of Technology,
Beijing 100124, China
(\email{xuzhenmath@bjut.edu.cn}).}
\and
Jianyuan Yin\footnotemark[2]
}

\usepackage{amsopn}
\DeclareMathOperator{\diag}{diag}

\ifpdf
\hypersetup{
  pdftitle={Trust-region high-index saddle dynamics},
  pdfauthor={G. Ji, L. Li, Z. Xu, and J. Yin}
}
\fi
\usepackage{amssymb,subfigure,float,verbatim,booktabs,makecell}

\usepackage{amsmath,enumitem,multirow}
\usepackage{algorithm,algpseudocode}
\usepackage[utf8]{inputenc}
\usepackage{textgreek}

\newtheorem{assumption}{Assumption}

\newcommand{\R}{\mathbb{R}}
\newcommand{\norm}[1]{\left\lVert #1 \right\rVert}

\begin{document}
\maketitle

\begin{abstract}
We propose a trust-region high-index saddle dynamics (TRHiSD) framework for computing high-index saddle points in ill-conditioned problems.
To address the deterioration of high-index saddle dynamics caused by inappropriate step sizes and search directions, we exploit its intrinsic local min--max structure and construct separate trust-region subproblems on the maximization and minimization subspaces to adapt the two trust-region radii to different local curvature scales. Under suitable assumptions and with exact spectral subspaces, we prove local linear convergence of the Cauchy-point variant and local quadratic convergence of the Newton-step variant. Numerical experiments on modified Rosenbrock functions are consistent with the predicted convergence rates. For the Landau--de Gennes model, a TRHiSD variant using a truncated conjugate-gradient step for the minimization subproblem substantially reduces the outer-iteration count and the wall-clock time in the tested cases.
\end{abstract}



\section{Introduction}

The computation of saddle points in complex systems has attracted considerable attention across many scientific disciplines. In physics, transition states are commonly characterized as index-1 saddle points of an energy functional and determine transition pathways between (meta)stable states. 
More generally, higher-index saddle points can also help reveal the global structure of the energy landscape and provide important information about bifurcations, competing transition mechanisms, and hierarchical connections among stationary states. 
Saddle-point searching methods have been successfully applied to identify excited states of Bose--Einstein condensates \cite{YIN2024100546,GAD_BEC} and to uncover defect structures, metastable configurations, and transition pathways in nematic liquid crystals described by the Landau--de Gennes (LdG) theory \cite{PhysRevLett,smectic_LDG}. 
They have also been used to study nucleation mechanisms and pattern formation in phase-field models \cite{math_pf,diblock_SL,triblock_SL,lattice_diffsion}.

A broad range of numerical methods have been developed for saddle-point computation. 
These methods can be generally classified into path-finding methods and surface-walking methods. 
Path-finding methods require both an initial state and a final state and seek a minimum energy path (MEP) connecting the two endpoints. 
When a converged MEP contains a single nondegenerate transition state, its maximum-energy configuration should correspond to an index-1 saddle point. 
Representative approaches include the string method  \cite{2002String,String_optimization} and the nudged elastic band method  \cite{NEB_Jonsson,NEB_Henkelman}. 
In contrast, surface-walking methods start from a single initial configuration and directly search for a nearby saddle point without requiring prior knowledge of the final state. 
These methods use local derivative information, such as gradients, Hessian-vector products, or approximations to relevant Hessian eigendirections. 
Representative examples include gentlest ascent methods and dynamics \cite{Crippen1971462,E_2011}, dimer-type methods \cite{dimer}, local minimax methods \cite{minimax_method,minimax_singularly}, and eigenvector-following methods \cite{CerjanMiller1981,Simons1983}, the iterative minimization formulation \cite{gao2015iterative},  gradient extremals \cite{jorgensen1988gradient} and Newton trajectories \cite{quapp1998searching}. 
This ability to search for saddle points without prescribing a final state makes surface-walking methods particularly suitable for exploring complex energy landscapes.

Among such methods, Yin et al. introduced the high-index saddle dynamics (HiSD) for computing index-$k$ saddle points \cite{2019High}. The discretization error and convergence properties of HiSD have subsequently been analyzed in \cite{HiSD_error,HiSD_Convergence}. Recently, Su et al. proposed an improved HiSD method and established its nonlocal convergence \cite{iHiSD}. The accelerated HiSD method introduced a momentum term to mitigate the slow convergence caused by ill-conditioning \cite{AHiSD}. Spectral step-size strategies, such as the Barzilai--Borwein (BB) method \cite{barzilai1988bb}, can partially account for local curvature information, but they do not provide an explicit mechanism for controlling the reliability of a trial step and may still produce excessively large steps in ill-conditioned regimes.

Therefore, while these approaches improve practical efficiency, they do not provide a model-based mechanism to determine whether a computed trial step is reliable. This limitation becomes particularly evident for ill-conditioned problems, where the Hessian spectrum contains widely separated curvature scales. Consequently, selecting a trial step that balances stability and efficiency becomes challenging. A conservative trial step to accommodate large curvature may lead to slow convergence, while aggressive trial steps may result in overshooting or unstable iterations.

Trust-region methods provide a well-established framework for controlling step sizes in nonlinear and nonconvex optimization \cite{trust_region_powell,trust_region_constrain,trust_region_free,trust_region_CG}. Instead of prescribing a step size directly, a trust-region method restricts a local model to a neighborhood in which the model is considered reliable and adjusts this neighborhood according to the agreement between the predicted and actual reductions.
This mechanism not only controls the step sizes but also helps reject unreliable trial steps while exploiting local curvature information, particularly in nonconvex problems. 
Moreover, the local min–max structure of HiSD allows the trust-region framework to be applied separately to the maximization and minimization subspaces. This enables independent adaptation of the two step components to their respective curvature scales and provides a natural alternative to explicit step-size selection in HiSD discretizations.
In particular, Newton-step and truncated conjugate-gradient (TCG) variants
use the block quadratic models to compute
curvature-informed search directions.

Near an index-\(k\) saddle point, HiSD exploits the local min--max structure by separating the dynamics into ascent and descent components associated with the unstable and stable spectral subspaces, respectively \cite{2019High}.
Motivated by this structure, we develop a trust-region HiSD (TRHiSD) method in which trust-region subproblems are solved separately in the two subspaces. 

The main contributions of this work are summarized as follows:
\begin{itemize}  
    \item We incorporate trust-region strategies into HiSD and develop a block-split trust-region framework with separate trust-region subproblems in the two subspaces. The trust-region mechanism is tailored to the saddle-point setting: the agreement between the predicted and actual increases is measured in the maximization subspace, whereas the agreement between the predicted and actual decreases is measured in the minimization subspace. The two trust-region radii are updated separately, allowing the ascent and descent components to adapt to different local curvature scales.

    \item Under suitable assumptions and with exact spectral subspaces, we establish a local convergence theory for the proposed framework. We derive a uniform lower bound for the trust-region radii and show that, sufficiently close to a nondegenerate index-\(k\) saddle point, the trust-region constraints become inactive and the first trial pair is accepted. This leads to local linear convergence of the Cauchy-point variant and local quadratic convergence of the Newton-step variant.

    \item For ill-conditioned problems, we further incorporate a TCG solver into the minimization trust-region subproblem to exploit curvature information effectively. Numerical experiments on modified Rosenbrock functions are consistent with the predicted local convergence behavior, while the LdG experiments show substantial reductions in outer iterations and computational time in the tested cases.

\end{itemize}

The remainder of this paper is organized as follows. \Cref{sec:2} reviews HiSD and presents the proposed trust-region framework together with several block-step strategies. \Cref{sec:3} establishes the local convergence theory for the Cauchy-point and Newton-step variants. \Cref{sec:4} reports numerical experiments on modified Rosenbrock functions and the LdG model. Conclusions and discussion for future work are presented in \Cref{sec:5}.

\section{TRHiSD framework}
\label{sec:2}
Let $F:\mathbb{R}^d\rightarrow\mathbb{R}$ be a twice continuously differentiable energy function, and define the gradient $g(x):=\nabla F(x)$, Hessian $H(x):=\nabla^2F(x).$ We use the Euclidean inner product $\langle u,v\rangle:=u^\top v$ in \(\R^d\) and the induced norm \(\norm{\cdot}\). For a matrix $A\in\mathbb{R}^{d\times d}$, we denote by $\|A\|:=\sup_{\|x\|=1}\|Ax\|$ the induced matrix 2-norm. Let \(x^*\) be a nondegenerate index-\(k\) saddle point of \(F\), that is, $g(x^*)=0$, and \(H(x^*)\) has exactly \(k\) negative eigenvalues and \(d-k\) positive eigenvalues according to Morse theory \cite{Milnor2016}.

\subsection{Review of HiSD}
HiSD is designed to locate index-\(k\) saddle points. 
The continuous HiSD system for locating an index-\(k\) saddle point of \(F\) is
\begin{equation}
\left\{
\begin{aligned}
\dot{x} &= -\left(I-2\sum_{i=1}^k v_i v_i^\top\right)g(x),\\
\dot{v}_i&=-\left(I-v_iv_i^\top-2\sum_{j=1}^{i-1}v_jv_j^\top\right)H(x)[v_i],\quad 1 \leq i\leq k.
\end{aligned}
\right.
\label{Hisd}
\end{equation}
Throughout this paper, \(I\) denotes the identity operator in the corresponding subspace. The continuous formulation is discussed in detail in \cite{2019High}, and the convergence rate and error estimates of the discrete algorithm are analyzed in \cite{HiSD_Convergence}.

The key idea underlying HiSD is to exploit the local min--max structure around an index-\(k\) saddle point $x^*$: $(v=0, w=0)$ is a solution to the min--max problem
\begin{equation}
    \label{exact-min-max-problem}
    \min\limits_{\substack{w\in W^*\\ \Vert w\Vert \leq \delta_1}}\max\limits_{\substack{v\in V^*\\ \Vert v \Vert \leq \delta_2}} F(x^*+w+v).
\end{equation}
Here, \(\delta_1>0\) and \(\delta_2>0\) are sufficiently small radii that restrict the minimization and maximization variables to local neighborhoods near the saddle point. Moreover, {$V^*$} is the $k$-dimensional subspace spanned by the eigenvectors of \(H(x^*)\) corresponding to its \(k\) negative eigenvalues, and \(W^*=(V^*)^\perp\) is the subspace associated with positive eigenvalues. This local min--max structure is reflected in the \(x\)-dynamics of HiSD: the gradient component in the subspace $V^*$ generates an ascent direction, while the component in the subspace $W^*$ generates a descent direction.

This observation provides the basic motivation for treating the local ascent and descent components separately when constructing the TRHiSD framework.

\subsection{Trust-region formulation}
Although fixed-step HiSD \cite{2019High} is straightforward to implement, its performance may deteriorate for ill-conditioned problems because a single step size cannot adapt to directions with substantially different local curvatures. Barzilai--Borwein (BB) step sizes can improve practical efficiency, but do not provide the same model-based step control as the trust-region methods and may require additional safeguards \cite{2019High,barzilai1988bb}. To improve robustness, we replace the explicit update of \(x^n\) by a trust-region method which automatically adjusts the steps. 

{At an iterate \(x^n\) sufficiently close to $x^*$, the eigenvalues of $H(x^n)$ can be ordered as $\lambda_1(x^n)\leq \cdots \lambda_k(x^n)<0<\lambda_{k+1}(x^n)\leq\cdots\leq \lambda_d(x^n)$ with the corresponding orthonormal eigenvectors $v_1(x^n),\cdots,v_d(x^n)$. For simplicity, we write $g^n:=g(x^n)$, $H^n:=H(x^n)$, $v_i^n:=v_i(x^n),\ 1\le i\le d$. Then the subspaces $V^n:=\operatorname{span}\{v_1^n,\ldots,v_k^n\}$ and $W^n:=(V^n)^\perp$ correspond locally to the ascent and descent components, respectively. We note that $V^n$ is spanned by eigenvectors corresponding to the $k$ smallest eigenvalues of $H^n$, with any orthonormal selection used in case of ties.} In a sufficiently small neighborhood of $x^*$, it coincides with the maximization subspace. 

Motivated by the local min--max structure in \eqref{exact-min-max-problem}, we consider
\begin{equation}
    \label{eq:local-min-max}
    \min_{p_2\in W^n}
    \max_{p_1\in V^n}
    F(x^n+p_1+p_2),
\end{equation}
where \(p_1\) and \(p_2\) are understood to be local perturbations. Directly solving the nonlinear local min--max problem is generally difficult. We therefore approximate \(F\) around \(x^n\) using the quadratic model
\begin{equation}
    \label{eq:model}
    m^n(p)
    :=
    F(x^n)+(g^n)^\top p+\frac12 p^\top H^n p .
\end{equation}
Introducing separate trust-region radii for the ascent and descent components leads to the local quadratic min--max problem, which can be decomposed into independent minimization and maximization subproblems:
\begin{equation}
    \label{eq:tr-min-max}
    \min_{\substack{p_2\in W^n\\
                    \|p_2\|\le\Delta_2^{n,\ell}}}
    \max_{\substack{p_1\in V^n\\
                    \|p_1\|\le\Delta_1^{n,\ell}}}
    m^n(p_1+p_2)=-F(x^n)+
    \min_{\substack{p_2\in W^n\\
                    \|p_2\|\le\Delta_2^{n,\ell}}} m^n(p_2)+
    \max_{\substack{p_1\in V^n\\
                    \|p_1\|\le\Delta_1^{n,\ell}}} m^n(p_1)      
\end{equation}
 where \(\Delta_1^{n,\ell}\) and \(\Delta_2^{n,\ell}\) are trust-region radii for the $\ell$-th inner trial at the fixed base point $x^n$. Since \(V^n\) and \(W^n\) are orthogonal invariant subspaces of \(H^n\), the cross term vanishes, so the quadratic model is separable with respect to \(p_1\) and \(p_2\). Accordingly, the quadratic min--max problem gives rise to the two block trust-region subproblems
\begin{equation}
\label{eq:block-tr-subproblems}
\min_{\substack{p_2\in W^n\\
                    \|p_2\|\le\Delta_2^{n,\ell}}} m^n(p_2),\quad
\max_{\substack{p_1\in V^n\\
                    \|p_1\|\le\Delta_1^{n,\ell}}} m^n(p_1).
\end{equation}
The trust-region strategies are then applied separately to the maximization and minimization subproblems to generate a pair of trial steps $p_1^{n,\ell}$ and $p_2^{n,\ell}$.

To assess the quality of the quadratic model in the two subspaces, we define the agreement ratios 
\begin{equation}
\label{eq:rho-inner}
\rho_1^{n,\ell}=\frac{F(x^n+p_1^{n,\ell})-F(x^n)}{\operatorname{Pred}_1^{n,\ell}},\qquad  \rho_2^{n,\ell}=\frac{F(x^n)-F(x^n+p_2^{n,\ell})}{\operatorname{Pred}_2^{n,\ell}},
\end{equation}
whenever the corresponding predicted increase or decrease is positive. Here, the predicted increase $\operatorname{Pred}_1^{n,\ell}:=m^n(p_1^{n,\ell})-m^n(0)$, and the predicted decrease $\operatorname{Pred}_2^{n,\ell}:=m^n(0)-m^n(p_2^{n,\ell})$. In the special case \(\operatorname{Pred}_i^{n,\ell}=0\), we set $\rho_i^{n,\ell}:=\frac12$.

If $\rho^{n,\ell}_1\geq 0$ and $\rho^{n,\ell}_2\geq0$, then both trial steps $p_1^{n,\ell}$ and $p_2^{n,\ell}$ are accepted. The next outer iterate is
\[
    x^{n+1}=x^n+p_1^{n}+p_2^{n},
\]
where we set $p_1^n=p_1^{n,\ell}$ and $p_2^n=p_2^{n,\ell}$. If at least one ratio $\rho_i^{n,\ell}$ is negative, then the trial pair is rejected and \(x^n\) remains the fixed base point. The inner iteration radii are updated as:
\begin{equation}
\label{eq:rejected-radius-update}
\Delta_i^{n,\ell+1}
=
\begin{cases}
    \Delta_i^{n,\ell},
    & \rho_i^{n,\ell}\geq 0,\\[1mm]
    \frac{1}{4}\Delta_i^{n,\ell},
    & \rho_i^{n,\ell}< 0,\\[1mm]
\end{cases}
\qquad i=1,2.
\end{equation}
This process is repeated until both agreement ratios are nonnegative. Once the trial pair is accepted, the radii associated with \(x^{n+1}\) are defined by
\begin{equation}
\label{eq:accepted-radius-update}
\Delta_i^{n+1}
=
\begin{cases}
    2\Delta_i^{n,\ell},
    & \rho_i^{n,\ell}>0.75,\\[1mm]
    \Delta_i^{n,\ell},
    & 0.25\le \rho_i^{n,\ell}\le0.75,\\[1mm]
    \frac14\Delta_i^{n,\ell},
    & 0\le \rho_i^{n,\ell}<0.25,
\end{cases}
\qquad i=1,2.
\end{equation}

Combining the block trust-region subproblems, the acceptance criterion, and the radius-update strategy, we obtain the TRHiSD framework summarized in \Cref{alg:trhisd}.

\begin{algorithm}[H]
\caption{TRHiSD framework for an index-\(k\) saddle point}
\label{alg:trhisd}
\begin{algorithmic}[1]
\Require
\(x^0\in\mathbb{R}^d\),
\(\Delta_1^0>0\),
\(\Delta_2^0>0\),
and outer iterations \(N\in\mathbb{N}\).
\For{\(n=0,1,\ldots,N-1\)}
    \State \label{step:subspace}Update the corresponding subspace \(V^n\) from \(H^n\) and set \(W^n=(V^n)^\perp\).
    \State Set $\Delta_i^{n,0}=\Delta_i^n$, $i=1,2$, and $\ell=0$
    \While{true}
        \State\label{step:block-step}Compute trial steps
        \(p_1^{n,\ell}\in V^n\) and
        \(p_2^{n,\ell}\in W^n\)
        based on subproblems
        \eqref{eq:block-tr-subproblems}.
        \State Compute the agreement ratios
        \(\rho_1^{n,\ell}\) and \(\rho_2^{n,\ell}\)
        from \eqref{eq:rho-inner}.
        \If{\(\rho_1^{n,\ell}\ge0\) and
             \(\rho_2^{n,\ell}\ge0\)}
            \State Accept the trial pair and set $x^{n+1}=x^n+p_1^{n,\ell}+p_2^{n,\ell}$.
            
            \State Update
            \(\Delta_i^{n+1}\), \(i=1,2\),
            according to rule \eqref{eq:accepted-radius-update}.
            \State \textbf{break}
        \Else
        \State Update \(\Delta_i^{n,\ell+1}\), \(i=1,2\), according to rule \eqref{eq:rejected-radius-update}.
        \State $\ell \leftarrow \ell+1$
        \EndIf
        
    \EndWhile
\EndFor
\State \Return \(x^N\).
\end{algorithmic}
\end{algorithm}

In Step \ref{step:subspace} of \Cref{alg:trhisd}, the subspace $V^n$ can be computed exactly or approximated by an iterative eigensolver initialized with the subspace from the previous iteration.
The convergence analysis in \Cref{sec:3} assumes that $V^n$ is an exact invariant subspace associated with the $k$ smallest eigenvalues of $H^n$, while numerical experiments use the locally optimal block preconditioned conjugate gradient (LOBPCG) method \cite{knyazev2001lobpcg} in \Cref{sec:4}.

\subsection{Step strategies}
The trial-step computation in Step~\ref{step:block-step} of \Cref{alg:trhisd} remains to be specified. Different choices of the block-step strategy lead to different TRHiSD variants. We consider three strategies: the Cauchy-point, Newton-step, and TCG-based inexact Newton-step variants.

To conveniently describe these block steps, let \(P^n:=\sum_{i=1}^k v_i^n (v_i^n)^\top\) and \(Q^n:=I-P^n\) denote the orthogonal projectors onto \(V^n\) and \(W^n\), respectively. Define $g_1^n:=P^n g^n$, $g_2^n:=Q^n g^n$, and the block Hessians $A_1^n:=-H^n|_{V^n}$, $A_2^n:=H^n|_{W^n}$. 
With this notation, we first define the Cauchy point.

\paragraph{Cauchy point}
The block Cauchy points are defined by
\[
    p_1^{n,\ell}
    =
    \alpha_1^{n,\ell} g_1^n,
    \qquad
    p_2^{n,\ell}
    =
    -\alpha_2^{n,\ell} g_2^n.
\]
If \(g_i^n\neq0\), the Cauchy step is
\begin{equation}
\label{eq:alpha_i}
\alpha_i^{n,\ell} =
\begin{cases}
\min\left\{
\frac{\|g_i^n\|^2}{(g_i^n)^\top A_i^n g_i^n},
\frac{\Delta_i^{n,\ell}}{\|g_i^n\|}
\right\},
& (g_i^n)^\top A_i^n g_i^n > 0,\\
\frac{\Delta_i^{n,\ell}}{\|g_i^n\|},
& (g_i^n)^\top A_i^n g_i^n \leq 0.
\end{cases}
\end{equation}
If \(g_i^n=0\), the corresponding Cauchy point is defined as
$$
p_i^{n,\ell} =
\begin{cases}
0,
&  A_i^n \succeq 0,\\
\Delta_i^{n,\ell} v_{\min}^n
& \text{else}.
\end{cases}
$$
where $v_{\min}^n$ is a unit eigenvector of $A_i^n$ associated with the smallest eigenvalue.

\paragraph{Newton step}
The Newton-step variant uses exact solutions of the trust-region subproblems \eqref{eq:block-tr-subproblems}. Since the feasible sets are compact and quadratic models are continuous, the exact solutions exist. When the block Hessians are positive definite and the unconstrained Newton steps lie inside the trust regions, these solutions coincide with the Newton steps.

\paragraph{Inexact Newton step}
When the block Newton equations admit approximate
solutions satisfying the prescribed residual tolerances,
the block inexact Newton steps satisfy
\[
A_1^n p_1^{n,\ell}-g_1^n=r_1^{n,\ell},
\qquad
A_2^n p_2^{n,\ell}+g_2^n=r_2^{n,\ell},
\]
subject to
\[
\|p_i^{n,\ell}\|\le\Delta_i^{n,\ell},
\qquad
\|r_i^{n,\ell}\|
\le\eta_i^{n,\ell}\|g_i^n\|,
\quad i=1,2.
\]
For general iterates, the block Hessians may be
indefinite or singular, so these residual conditions
are not guaranteed to be attainable.

A practical approach for computing a trust-region step is the TCG method \cite{trust_region_CG}. We apply TCG to the minimization trust-region subproblem
\begin{equation}
\label{eq:tcg-subproblem}
\min_{\substack{p\in W^n\\
\|p\|\le\Delta_2^{n,\ell}}}
\left\{
(g_2^n)^\top p+\frac12p^\top A_2^n p
\right\}.
\end{equation}
This subproblem admits a global minimizer even when
$A_2^n$ is indefinite or singular, since its feasible
set is compact and its objective is continuous. For $g_2^n=0$, we set $p_2^{n,\ell}=0$ if $A_2^n\succeq0$. Otherwise, choose a unit negative-curvature direction $p\in W^n$ satisfying $p^\top A_2^n p<0$ and set $p_2^{n,\ell}=\Delta_2^{n,\ell}p$. For $g_2^n\neq0$, we apply the TCG iteration. The TCG iteration starts from $p^{(0)}=0$ and terminates when any of the following conditions is met: the residual criterion
\[
\|A_2^np^{(j)}+g_2^n\|
\le\epsilon_{\rm tcg}\|g_2^n\|
\]
is satisfied, the trust-region boundary is reached,
the maximum number of TCG iterations $i_{\max}$
is attained, or nonpositive curvature is detected.
If the residual criterion is satisfied before the
trust-region boundary is reached and before nonpositive
curvature is detected, the resulting step is an
inexact Newton step in the above sense. Otherwise,
the computed step is treated as an approximate
trust-region step.

We denote by TRHiSD-CP the variant that uses Cauchy points on both \(V^n\) and \(W^n\), by TRHiSD-N the variant that uses Newton steps on both subspaces, and by TRHiSD-CP/TCG the hybrid variant that uses a Cauchy point on \(V^n\) and a TCG step on \(W^n\).

\begin{remark}
Fix an outer iterate \(x^n\). Suppose that \(F\) is twice continuously differentiable in a neighborhood of \(x^n\), and that the Cauchy-point and Newton-step variant satisfy $\operatorname{Pred}_i^{n,\ell}\ge c_n\|p_i^{n,\ell}\|^2$, $i=1,2$, for some \(c_n>0\). {This lower bound is established locally in \Cref{lem:pred-bound} under \Cref{ass:basic}.} Then the inner iteration terminates after finitely many trials. Indeed, $F(x^n+p)-m^n(p)=o(\|p\|^2)$, so whenever \(p_i^{n,\ell}\neq0\), the agreement ratios satisfy $\rho_i^{n,\ell}\to1$\ (as $\Delta_i^{n,\ell}\to0$). Since every rejected block has its radius reduced by a factor \(1/4\), each block can be rejected only a finite number of times. If \(p_i^{n,\ell}=0\), then \(\rho_i^{n,\ell}=1/2\) by convention and the block is not rejected. 
\end{remark}

\section{Local convergence analysis}
\label{sec:3}
In this section, we establish the local convergence of TRHiSD. After introducing the assumptions and auxiliary estimates, we derive a uniform lower bound for the trust-region radii before the first exit from a local neighborhood and show that, sufficiently close to a nondegenerate index-$k$ saddle point, the trust-region constraints become inactive and the first trial pair is accepted. We then establish local linear convergence of the Cauchy-point variant and local quadratic convergence of the Newton-step variant.

\subsection{Notation and assumptions}
We retain the notation introduced in \Cref{sec:2} and additionally denote $e^n:=x^n-x^*$. To prove the local convergence, we introduce the following assumptions, which are standard in the local convergence analysis of iterative methods.
\begin{assumption}[Local smoothness and spectral gap]\label{ass:basic}
There exist a convex neighborhood \(U\) of the index-$k$ saddle point \(x^*\), constants \(M>0\) and
\(L\ge \mu>0\), such that the following statements hold.
\begin{enumerate}[label=(A\arabic*)]
\item \(F\in C^2(U)\), and the Hessian is Lipschitz continuous on \(U\):
\[
    \norm{H(x)-H(y)}\le M\norm{x-y},
    \qquad x,y\in U.
\]
\item For every \(x\in U\), \(H(x)\) has exactly \(k\) negative eigenvalues
and \(d-k\) positive eigenvalues. Moreover, the following bounds hold:
\[
    \lambda_i(H(x))\le -\mu,
    \; 1\le i\le k;
    \quad
    \lambda_i(H(x))\ge \mu,
    \; k+1\le i\le d;
    \quad
    \norm{H(x)}\le L.
\]
\end{enumerate}
\end{assumption}

If $x^n\in U$, by \Cref{ass:basic}, the block Hessians $A_1^n$ and $A_2^n$ are positive definite and satisfy $\mu I\preceq A_i^n\preceq LI$ for $i=1,2$. To prove the local convergence of TRHiSD, we introduce the absolute Hessian
\[
    B^n:=-P^nH^nP^n+Q^nH^nQ^n.
\]
Under \Cref{ass:basic}, $B^n$ is positive definite and satisfies $\mu I\preceq B^n\preceq LI$.
We can define the \(B^n\)-norm as $\norm{z}_{B^n}:=(z^\top B^n z)^{1/2}$ which satisfies 
$
    \sqrt{\mu}\norm{z}
    \le
    \norm{z}_{B^n}
    \le
    \sqrt{L}\norm{z}
$.

\subsection{Auxiliary lemmas}
We first collect three standard estimates used in the local convergence
analysis: a Taylor remainder bound, the exact-line-search steepest-descent contraction, and the local Lipschitz continuity of the absolute Hessian. Since $g(x^*)=0$, we can express $g^n$ in terms of $H^n e^n$ and a Taylor remainder. 
Define the Taylor remainder

\begin{equation}
    r^n:=g^n-H^n e^n.
\end{equation}
The following lemma provides a quadratic bound for this remainder and establishes the corresponding
relations for the projected gradient components.

\begin{lemma}[Taylor remainder]\label{lem:taylor}
Under \Cref{ass:basic}, for \(x^n\in U\), we have
\[
    \norm{r^n}\le \frac{M}{2}\norm{e^n}^2.
\]
Moreover, if $z_1^n:=P^n e^n$, $z_2^n:=Q^n e^n$,
then we have 
\[
    g_1^n=-A_1^n z_1^n+P^n r^n,
    \qquad
    g_2^n=A_2^n z_2^n+Q^n r^n.
\]
\end{lemma}

\begin{proof}
Since \(g(x^*)=0\), $g^n=g^n-g^*=\int_0^1 H(x^*+t e^n)e^n\,dt$. Therefore, we define
\[
    r^n
    :=
    g^n-H^n e^n
    =
    \int_0^1\bigl(H(x^*+t e^n)-H^n\bigr)e^n\,dt .
\]
The Lipschitz continuity of \(H\) gives
\[
    \norm{r^n}
    \le
    \int_0^1 M(1-t)\norm{e^n}^2\,dt
    =
    \frac{M}{2}\norm{e^n}^2.
\]
The projected identities follow from the invariance of \(V^n\) and \(W^n\) under \(H^n\).
\end{proof}

To analyze the convergence, the following lemmas reveal the recursion relationship of $x^{n+1}-x^*$ and $x^n-x^*$.
\begin{lemma}[Steepest-descent contraction]\label{lem:sd}
Let \(A\) be symmetric positive definite and satisfy
$\mu I\preceq A\preceq LI$. 
For \(w\neq0\), define
\[
    \alpha(w):=
    \frac{\norm{Aw}^2}{(Aw)^\top A(Aw)},
    \qquad
    w^+:=w-\alpha(w)Aw.
\]
Then we have
\[
    \norm{w^+}_A\le q\norm{w}_A,
    \qquad
    q:=\frac{L-\mu}{L+\mu}<1,
\]
where \(\norm{z}_A:=(z^\top A z)^{1/2}\).
\end{lemma}

\begin{proof}
This is the classical contraction estimate for the steepest-descent method with exact line search applied to a strongly convex quadratic function; see, for example, \cite[Theorem~3.3]{NocedalWright2006}.
\end{proof}

\begin{lemma}[Local Lipschitz continuity of the absolute Hessian]
\label{lem:B-lipschitz}
Under \Cref{ass:basic}, the absolute Hessian
\(B(x)=\sum_{i=1}^d |\lambda(x)|v_i(x)v_i(x)^\top\) is locally Lipschitz continuous on \(U\). More precisely, there exists a constant
\(L_B>0\) such that
\[
    \|B(x)-B(y)\|
    \le
    L_B\|x-y\|,
    \qquad x,y\in U.
\]
\end{lemma}
\begin{proof}
Let \(P(x):=\sum_{i=1}^k v_i(x)v_i(x)^\top\) be the spectral projector of $H(x)$ associated with its negative eigenvalues. Choose a positively oriented closed contour \(\Gamma\) in the complex plane which encloses the interval \([-L,-\mu]\) and does not intersect the spectrum
of \(H(x)\) for any \(x\in U\). Since the spectral gap is uniform, there exists a constant \(d_\Gamma>0\) such that $\operatorname{dist}\bigl(z,\sigma(H(x))\bigr)\ge d_\Gamma,\ z\in\Gamma,\ x\in U$. It follows that
\[
    \|(zI-H(x))^{-1}\|
    \le
    \frac1{d_\Gamma},
    \qquad
    z\in\Gamma,\quad x\in U.
\]

By the Riesz projection formula and the resolvent identity (see, e.g., \cite{Kato1995}), we obtain
\[
\begin{aligned}
    \|P(x)-P(y)\|
    &\le
    \frac{1}{2\pi}
    \int_\Gamma
    \|(zI-H(x))^{-1}\|
    \|H(x)-H(y)\|
    \|(zI-H(y))^{-1}\|
    |dz|  \\
    &\le
    \frac{|\Gamma|}{2\pi}
    \frac{1}{d_\Gamma^2}
    \|H(x)-H(y)\| .
\end{aligned}
\]
Since \(H\) is \(M\)-Lipschitz on \(U\),  we obtain $\|P(x)-P(y)\|\leq C_P\|x-y\|$, where $C_P:=\frac{|\Gamma|}{2\pi}\frac{M}{d_\Gamma^2}$.

Since
$
    B(x)
    =
    -P(x)H(x)P(x)+Q(x)H(x)Q(x)
    =(I-2P(x))H(x),
$
we write
\[
\begin{aligned}
    B(x)-B(y)
    &=
    (I-2P(x))H(x)
    -
    (I-2P(y))H(y) \\
    &=
    (I-2P(x))(H(x)-H(y))
    -
    2(P(x)-P(y))H(y).
\end{aligned}
\]
Taking norms gives
\[
\begin{aligned}
    \|B(x)-B(y)\|
    &\le
    \|I-2P(x)\|\,\|H(x)-H(y)\|
    +
    2\|P(x)-P(y)\|\,\|H(y)\|.
\end{aligned}
\]
Since \(P(x)\) is an orthogonal projector, \(I-2P(x)\) is an orthogonal reflection and therefore satisfies $\|I-2P(x)\|=1$. Since $\|H(y)\|\le L$ by \Cref{ass:basic}, we obtain
\[
    \|B(x)-B(y)\|
    \le
    L_B\|x-y\|,
\]
where $L_B:=M+2LC_P$.
This proves that $B(x)$ is locally Lipschitz on \(U\).
\end{proof}

\subsection{Lower bound for the radii}
Fix \(R>0\) such that $\overline{B}(x^*,R)\subset U$. We define the first exit time $\tau$ from \(B(x^*,R/2)\) by
\[
    \tau
    :=
    \inf\{n\ge0:\ x^n\notin B(x^*,R/2)\}.
\]
All statements below that involve radii are first proved for \(n<\tau\). The
bootstrap argument in the convergence theorems then shows that \(\tau=\infty\) whenever the initial iterate $x^0$ is sufficiently close to $x^*$.
\begin{lemma}[Block predicted-variation bound]\label{lem:pred-bound}
Let \(x^n\in U\). For either the Cauchy point or Newton step, the block predicted variations satisfy
\[
    \operatorname{Pred}_i^{n,\ell}\ge \frac{\mu}{2}\norm{p_i^{n,\ell}}^2,
    \qquad i=1,2.
\]
\end{lemma}

\begin{proof}
For clarity, we present the proof for the maximization subproblem ($i=1$), and the case $i=2$ follows analogously. For the Cauchy-point variant, we write
\(p_1^{n,\ell}=\alpha g_1^n\), where $0\le \alpha\le\frac{\norm{g_1^n}^2}{(g_1^n)^\top A_1^n g_1^n}\le \frac1\mu$. Then, we obtain
\[
    \operatorname{Pred}_1^{n,\ell}
    =
    \alpha\norm{g_1^n}^2
    -\frac12\alpha^2 (g_1^n)^\top A_1^n g_1^n  
    \ge
    \frac12\alpha\norm{g_1^n}^2
    \ge
    \frac{\mu}{2}\alpha^2\norm{g_1^n}^2
    =
    \frac{\mu}{2}\norm{p_1^{n,\ell}}^2 .
\]
For the exact Newton step, the result follows from the standard optimality
conditions. For the maximization subproblem ($i=1$), the model on \(V^n\) is
\[
    m^n(p)-m^n(0)
    =
    (g_1^n)^\top p-\frac12 p^\top A_1^n p.
\]
The optimality conditions give
\(
    (A_1^n+\sigma I)p_1^{n,\ell}=g_1^n,
    \ \sigma\ge0.
\)
Therefore, we obtain
\[
    \operatorname{Pred}_1^{n,\ell}
    =
    \frac12(p_1^{n,\ell})^\top A_1^n p_1^{n,\ell}
    +\sigma\norm{p_1^{n,\ell}}^2
    \ge
    \frac{\mu}{2}\norm{p_1^{n,\ell}}^2.
\]
\end{proof}

\begin{lemma}[Uniform lower bound for the radii before exit]
\label{lem:radius-lower}
Under \Cref{ass:basic}, 
let $\delta_{\rm sh}
    :=
    \min\left\{
        \frac{9\mu}{4M},
        \frac{R}{2}
    \right\}$.
For arbitrary prescribed initial radii
$\Delta_1^0>0$, $\Delta_2^0>0$,
define
\[
    \underline\Delta_i
    :=
    \min\left\{
        \Delta_i^0,
        \frac14\delta_{\rm sh}
    \right\},
    \qquad i=1,2.
\]
Then, for every outer iteration \(n<\tau\) and every inner trial \(\ell\)
generated before accepting \(x^{n+1}\), we have 
$\Delta_i^{n,\ell}\ge \underline\Delta_i$ for $i=1,2$, and consequently, $\Delta_i^n\ge \underline\Delta_i$ for $i=1,2$.
\end{lemma}

\begin{proof}
Fix a block \(i\). Suppose \(n<\tau\), so \(x^n\in B(x^*,R/2)\). If $\Delta_i^{n,\ell}\le \delta_{\rm sh}$, 
then the trial step satisfies
\(
    \norm{p_i^{n,\ell}}
    \le
    \Delta_i^{n,\ell}
    \le
    \delta_{\rm sh}
    \le
    \frac{R}{2}.
\)
Therefore the whole segment \(x^n+t p_i^{n,\ell}\), \(0\le t\le1\), lies in
\(B(x^*,R)\). Taylor's theorem and the Lipschitz continuity of \(H\) give
\[
    \left|
        F(x^n+p_i^{n,\ell})-m^n(p_i^{n,\ell})
    \right|
    \le
    \frac{M}{6}\norm{p_i^{n,\ell}}^3.
\]
By \Cref{lem:pred-bound}, we have $\operatorname{Pred}_i^{n,\ell}
    \ge
    \frac{\mu}{2}\norm{p_i^{n,\ell}}^2$. For $\operatorname{Pred}_i^{n,\ell}>0$, we have,
\[
    |\rho_i^{n,\ell}-1|
    \le
    \frac{(M/6)\norm{p_i^{n,\ell}}^3}
         {(\mu/2)\norm{p_i^{n,\ell}}^2}
    =
    \frac{M}{3\mu}\norm{p_i^{n,\ell}}
    \le
    \frac{M}{3\mu}\Delta_i^{n,\ell}.
\]
Hence, if \(\Delta_i^{n,\ell}\le\delta_{\rm sh}\), then
\(
    \rho_i^{n,\ell}
    \ge
    1-\frac{M}{3\mu}\delta_{\rm sh}
    \ge
    \frac14.
\)
For \(\operatorname{Pred}_i^{n,\ell}=0\), the same conclusion holds by convention \(\rho_i^{n,\ell}=1/2\).

Therefore, whenever the radius of block \(i\) is at most
\(\delta_{\rm sh}\), it is neither rejected nor reduced. Consequently, any future reduction of this radius can occur only from a value strictly larger than \(\delta_{\rm sh}\). Such a reduction multiplies the radius by \(1/4\), and hence the new radius is strictly larger than \(\delta_{\rm sh}/4\). Starting from \(\Delta_i^0\), induction over all outer and inner updates before the exit time gives $\Delta_i^{n,\ell}\ge\min\left\{\Delta_i^0,\frac14\delta_{\rm sh}\right\}=\underline\Delta_i$.
\end{proof}

The preceding lemma shows that the trust-region radii have a lower bound. We next show that, sufficiently close to \(x^*\), this lower bound makes the trust-region constraints inactive and guarantees first-trial acceptance.

\subsection{Local untruncated-step region}
Since the initial trust-region radii \(\Delta_1^0\) and \(\Delta_2^0\) are positive, \Cref{lem:radius-lower} ensures that $\underline\Delta>0$. We define
\begin{equation}
    r_{\rm tr}:=\frac{\mu}{2L}\underline\Delta.
\end{equation}
The following lemma shows that, within $B(x^*,r_{\rm tr})$, the trust-region constraints become inactive and the first trial pair is accepted. Under \Cref{ass:basic}, the block Hessians $A_i^n$ are positive definite for $i=1,2$. When the trust-region constraints are inactive, the Cauchy steps in \eqref{eq:alpha_i} reduce to the untruncated Cauchy steps
\begin{equation}
\label{step:untruncated_C}
 p_1^{C,n}
=\alpha_1^n g_1^n,\quad
\alpha_1^n
=\frac{\|g_1^n\|^2}
{(g_1^n)^\top A_1^n g_1^n},\quad
p_2^{C,n}
=-\alpha_2^n g_2^n,\quad
\alpha_2^n=
\frac{\|g_2^n\|^2}
{(g_2^n)^\top A_2^n g_2^n},
\end{equation}
Similarly, the corresponding untruncated Newton steps are given by
\begin{equation}
\label{step:untruncated_N}
    p_1^{N,n}:=(A_1^n)^{-1}g_1^n,
    \quad
    p_2^{N,n}:=-(A_2^n)^{-1}g_2^n.
\end{equation}

\begin{lemma}[Local untruncated-step region and first-trial acceptance]
\label{lem:pure-local}
Under \Cref{ass:basic}, if \(n<\tau\) and $x^n\in B(x^*,r_{\rm tr})$,  then both the Cauchy point and the split Newton step are strictly inside their corresponding trust-region balls. In this case, the trust-region constraints are inactive and the untruncated Cauchy point and Newton step are accepted.
\end{lemma}

\begin{proof}
By \Cref{lem:radius-lower}, we have $\Delta_i^{n}\ge \underline\Delta$ for $i=1,2$.
Since \(g(x^*)=0\) and \(\norm{H(x)}\le L\) in \(U\), we have 
\[
    \norm{g^n}
    =
    \norm{g(x^n)-g(x^*)}
    \le
    L\norm{x^n-x^*}.
\]

For the untruncated Cauchy point, we have
\[
    \norm{p_i^{C,n}}
    \le
    \frac1\mu\norm{g_i^n}
    \le
    \frac{L}{\mu}\norm{x^n-x^*},
    \qquad i=1,2.
\]
Hence the same argument shows \(\norm{p_i^{C,n}}<\Delta_i^{n}\). Thus the Cauchy point is also strictly interior. For the untruncated Newton step, we have
\[
    \norm{p_i^{N,n}}
    \le
    \frac{L}{\mu}\norm{x^n-x^*},
    \qquad i=1,2.
\]
If \(x^n\in B(x^*,r_{\rm tr})\), then $\norm{p_i^{N,n}}<\frac{L}{\mu}r_{\rm tr}=
\frac12\underline\Delta\le\Delta_i^{n}$. Thus the split Newton step is strictly interior.

Let \(p_i^{n,0}\) denote either the Cauchy point or Newton step. From the preceding estimates, we obtain $\|p_i^{n,0}\|<\frac{\underline{\Delta}_i}{2}\leq\frac{\delta_{\mathrm{sh}}}{8}\leq\frac{R}{16}$. Moreover, we have $\|x^n-x^*\|<r_{\mathrm{tr}}\leq\frac{\Delta}{2}\leq\frac{R}{16}.$
Hence, we have
\(
    x^n+t p_i^{n,0}\in B(x^*,R),
    \ 0\leq t\leq1.
\)
Taylor's theorem and \Cref{lem:pred-bound} therefore give
\[
    |\rho_i^{n,0}-1|
    \leq
    \frac{(M/6)\|p_i^{n,0}\|^3}
         {(\mu/2)\|p_i^{n,0}\|^2}  
    =
    \frac{M}{3\mu}\|p_i^{n,0}\|  
    \leq
    \frac{M}{6\mu}\Delta
    \leq
    \frac{M}{24\mu}\delta_{\mathrm{sh}}
    \leq
    \frac{3}{32}.
\]
whenever \(p_i^{n,0}\neq0\). Thus, we obtain $\rho_i^{n,0}\geq\frac{29}{32}>\frac34$. If \(p_i^{n,0}=0\), then \(\rho_i^{n,0}=1/2\) by convention.
Consequently, both block ratios are nonnegative and the first trial pair is
accepted.
\end{proof}

\Cref{lem:pure-local} reduces the local TRHiSD iterations to the corresponding
untruncated Cauchy and Newton iterations. We now analyze these two cases separately.

\subsection{Local linear convergence of TRHiSD-CP}
By \Cref{lem:pure-local}, sufficiently close to the target saddle point, the first trial pair is accepted and the Cauchy steps coincide with the untruncated steps defined in \eqref{step:untruncated_C}. Consequently, the TRHiSD-CP iteration reduces locally to
$
    x^{n+1}=x^n+p_1^{C,n}+p_2^{C,n}.
$
\begin{lemma}[Bounds on the untruncated Cauchy steps]\label{lem:cauchy-step}
Every untruncated Cauchy coefficient $\alpha_i^n = \frac{\Vert g_i^n \Vert^2}{g_i^n A_i^n g_i^n}$ satisfies
\(
    \frac1L\le \alpha_i^n\le \frac1\mu,\ i=1,2,
\)
whenever \(g_i^n\neq0\). Moreover, the step satisfies $\norm{x^{n+1}-x^n}\le\frac{L}{\mu}\norm{e^n}$.
\end{lemma}

\begin{proof}
The bounds on \(\alpha_i^n\) follow from the Rayleigh quotient bounds
\[
    \mu\norm{v}^2\le v^\top A_i^n v\le L\norm{v}^2.
\]
Moreover, we have
\(
    \norm{g^n}
    =
    \norm{g(x^n)-g(x^*)}
    \le
    L\norm{e^n}.
\)
Since \(p_1^{C,n}\perp p_2^{C,n}\), it follows that
\[
    \norm{x^{n+1}-x^n}^2
    =
    \norm{p_1^{C,n}+p_2^{C,n}}^2
    \le
    \frac1{\mu^2}\left(\norm{g_1^n}^2+\norm{g_2^n}^2\right)
    =
    \frac1{\mu^2}\norm{g^n}^2
    \le
    \frac{L^2}{\mu^2}\norm{e^n}^2.
\]
\end{proof}
\begin{theorem}[Local linear convergence]
\label{thm:cauchy}
Let \(\Delta_1^0>0\), \(\Delta_2^0>0\) and \(\eta\in(q,1)\) be arbitrary, where $q=\frac{L-\mu}{L+\mu}$. Under
\Cref{ass:basic}, there exists
\(0<\delta_C<\min\{\sqrt{\mu}\frac{R}{2},\sqrt{\mu} r_{tr}\}\), depending on \(R,\Delta_1^0,\Delta_2^0,M,L,\mu,L_B,\eta\), such that, if $\norm{x^0-x^*}_{B^0}\le \delta_C$, then the TRHiSD-CP iteration is well defined, all accepted steps are the untruncated Cauchy points, and
\[
    \norm{x^{n+1}-x^*}_{B^{n+1}}
    \le
    \eta\norm{x^n-x^*}_{B^n},
    \qquad n\geq 0.
\]
Thus \(x^n\to x^*\) linearly.
\end{theorem}
\begin{proof}
Set $y_n:=\|e^n\|_{B^n}$, and define
\(
    \bar{\delta}
    :=
    \min
    \left\{
        \frac{\sqrt{\mu}R}{2},
        \sqrt{\mu}r_{\mathrm{tr}}
    \right\}.
\)
Suppose that \(y_n\leq \bar{\delta}\). Then, we have $\|e^n\|\leq\frac{y_n}{\sqrt{\mu}}\leq\min\left\{\frac{R}{2},r_{\mathrm{tr}}\right\}.$ \Cref{lem:pure-local} implies that the first trial pair is accepted and that the accepted Cauchy points are untruncated. Since $e^n=z_1^n+z_2^n$, the iteration gives
\[
e^{n+1}=e^n+\alpha_1^ng_1^n-\alpha_2^ng_2^n=(I-\alpha_1^nA_1^n)z_1^n
 +(I-\alpha_2^nA_2^n)z_2^n+\alpha_1^nP^nr^n-\alpha_2^nQ^nr^n.
\]
To separate the linear contraction terms from the
Taylor remainder, we define $w_1^n:=z_1^n-(A_1^n)^{-1}P^nr^n$,and $w_2^n:=z_2^n+(A_2^n)^{-1}Q^nr^n$. Then, we have
\[
e^{n+1}
=(I-\alpha_1^nA_1^n)w_1^n+(I-\alpha_2^nA_2^n)w_2^n+(A_1^n)^{-1}P^nr^n-(A_2^n)^{-1}Q^nr^n.
\]
By \Cref{lem:sd}, we have $\|(I-\alpha_i^nA_i^n)w_i^n\|_{A_i^n} \leq q\|w_i^n\|_{A_i^n},\ i=1,2$. Moreover, we have $\|(A_1^n)^{-1}P^nr^n\|_{A_1^n}^2+\Vert(A_2^n)^{-1}Q^nr^n\|_{A_2^n}^2=\|r^n\|_{(B^n)^{-1}}^2$. Therefore, by the triangle inequality, we have
\[
\|e^{n+1}\|_{B^n}\leq{}q\left(\|w_1^n\|_{A_1^n}^2+\|w_2^n\|_{A_2^n}^2\right)^{1/2}+\|r^n\|_{(B^n)^{-1}}.
\]
Moreover, we have the estimate
\[
    \left(
        \|w_1^n\|_{A_1^n}^2
        +
        \|w_2^n\|_{A_2^n}^2
    \right)^{1/2}
    \leq
    y_n+\|r^n\|_{(B^n)^{-1}}.
\]
Hence, we obtain
\[
    \|e^{n+1}\|_{B^n}
    \leq
    qy_n+(1+q)\|r^n\|_{(B^n)^{-1}}
    \leq
    qy_n+C_0y_n^2,
\]
where $C_0:=\frac{(1+q)M}{2\mu^{3/2}}$.

We note that \Cref{lem:cauchy-step} and
\(\|e^n\|\leq r_{\mathrm{tr}}\) yield
\(
    \|x^{n+1}-x^n\|
    \leq
    \frac12\Delta
    \leq
    \frac{R}{16}.
\)
Since \(\|x^n-x^*\|\leq R/2\), we have
\[
    \|x^{n+1}-x^*\|
    \leq
    \|x^n-x^*\|
    +
    \|x^{n+1}-x^n\| 
    \leq
    \frac{R}{2}+\frac{R}{16}
    <R.
\]
Thus \(x^{n+1}\in B(x^*,R)\), and \Cref{lem:B-lipschitz} is applicable. Furthermore, the Lipschitz continuity of $B$ gives
\[
    \|B^{n+1}-B^n\|
    \leq
    L_B\|x^{n+1}-x^n\|
    \leq
    \frac{L_BL}{\mu^{3/2}}y_n .
\]

Moreover, since $B^n\succeq\mu I$, we obtain
\[
\begin{aligned}
\|e^{n+1}\|_{B^{n+1}}^2
&\leq \|e^{n+1}\|_{B^n}^2
+\|B^{n+1}-B^n\|_2\|e^{n+1}\|_2^2\\
&\leq
\left(1+\frac{\|B^{n+1}-B^n\|_2}{\mu}\right)
\|e^{n+1}\|_{B^n}^2\\
&\leq
(1+K_By_n)\|e^{n+1}\|_{B^n}^2.
\end{aligned}
\]
Taking square roots and using $\sqrt{1+t}\leq 1+t$ for $t\geq0$, we obtain the desired estimate.
Hence, we have
\[
    \|e^{n+1}\|_{B^{n+1}}
    \leq
    (1+K_By_n)
    \|e^{n+1}\|_{B^n},
\]
where $K_B:=\frac{L_BL}{\mu^{5/2}}$.

Combining the preceding estimates, we obtain
\begin{equation}
    y_{n+1}
    \leq
    (1+K_By_n)
    (qy_n+C_0y_n^2)
    =
    qy_n
    +(C_0+K_Bq)y_n^2
    +K_BC_0y_n^3 .
\label{eq:cauchy-recursion}
\end{equation}
Define \(C_1:=C_0+K_Bq, \ C_2:=K_BC_0\). Then we have
\begin{equation}
\label{eq:cauchy-recursion-simple}
    y_{n+1}
    \leq
    qy_n+C_1y_n^2+C_2y_n^3 .
\end{equation}
Choose $\delta_C\leq \bar{\delta}$ such that
\(
    C_1\delta_C+C_2\delta_C^2
    \leq
    \eta-q .
\)

We now use a bootstrap induction. 
The initial condition gives
$\|e^0\|\leq y_0/\sqrt{\mu}<R/2$. Supposing {$n<\tau$ and} $y_n\leq\delta_C$, all preceding estimates remain valid, and
\eqref{eq:cauchy-recursion-simple} gives
\[
    y_{n+1}
    \leq
    (q+C_1\delta_C+C_2\delta_C^2)y_n
    \leq
    \eta y_n
    \leq
    \delta_C .
\]
Since $x^{n+1}\in B(x^*,R)$, norm equivalence
yields $\|e^{n+1}\|\leq y_{n+1}/\sqrt{\mu}<R/2$.
Thus $n+1<\tau$, and the induction closes.
Starting from \(y_0\leq\delta_C\), induction yields $y_n\leq\eta^n y_0,\ n\geq0$. Therefore, $x^n\rightarrow x^*$ linearly.
\end{proof}

\subsection{Local quadratic convergence of TRHiSD-N}
By \Cref{lem:pure-local}, sufficiently close to the target saddle point, the first trial pair is accepted and the Newton steps coincide with the
untruncated Newton steps defined in \eqref{step:untruncated_N}. Consequently, the TRHiSD-N iteration locally reduces to $x^{n+1}=x^n+p_1^{N,n}+p_2^{N,n}$.

\begin{lemma}[Equivalence of split and full Newton steps]\label{lem:split-newton}
The untruncated Newton steps in \eqref{step:untruncated_N} satisfy
$
    p_1^{N,n}+p_2^{N,n}=-(H^n)^{-1}g^n.
$
\end{lemma}
\begin{proof}

On \(V^n\), \(H^n=-A_1^n\). Hence, we have $p_1^{N,n}=(A_1^n)^{-1}g_1^n=-(H^n)^{-1}P^n g^n$. The analogous identity on $W^n$ follows similarly. Since $P^n+Q^n=I$, adding the two identities proves the result.
\end{proof}

\begin{theorem}[Local quadratic convergence]
\label{thm:newton}
Let \(\Delta_1^0>0\) and \(\Delta_2^0>0\) be arbitrary. Under
\Cref{ass:basic}, there exists \(\delta_N=\min\left\{\frac{R}{2},r_{\rm tr},\frac{\mu}{M}\right\}>0\), depending on
\(\Delta_1^0,\Delta_2^0,M,L,\mu,R\), such that if \(\norm{x^0-x^*}< \delta_N\), then the Newton-step TRHiSD iteration is well defined, all accepted steps are the untruncated split Newton steps, and
\[
    \norm{x^{n+1}-x^*}
    \le
    \frac{M}{2\mu}\norm{x^n-x^*}^2,
    \qquad n\geq 0.
\]
Consequently, \(x^n\to x^*\) quadratically.
\end{theorem}

\begin{proof}
We prove \(\norm{e^n}< \delta_N\) for all \(n\ge0\) by induction.
This is true for \(n=0\). Suppose it is true at step \(n\). Then
\(x^n\in B(x^*,R/2)\), so \(n<\tau\). Since \(\delta_N\le r_{\rm tr}\),
\Cref{lem:pure-local} implies that the trust-region constraints are inactive and the steps are accepted.

By \Cref{lem:split-newton}, we have $e^{n+1}=e^n-(H^n)^{-1}g^n$. Using \Cref{lem:taylor}, we obtain
\[
    e^{n+1}
    =
    e^n-(H^n)^{-1}(H^n e^n+r^n)
    =
    -(H^n)^{-1}r^n.
\]
Since $\|(H^n)^{-1}\|\leq\frac{1}{\mu}$, we obtain
\[
    \norm{e^{n+1}}
    \le
    \frac1\mu\norm{r^n}
    \le
    \frac{M}{2\mu}\norm{e^n}^2
    \le
    \frac12\norm{e^n}
    <
    \delta_N.
\]
The induction closes. Hence, the iterates never exit the local ball $B(x^*,\frac{R}{2})$, the trust-region constraints remain inactive, and the quadratic estimate holds for every \(n\).
\end{proof}

\begin{remark}
\label{rem:inexact-eigenspace}
The preceding analysis assumes that $v_1^n,...,v_k^n$ are exact eigenvectors of $H(x^n)$. In practice, let \(\widetilde H^n\) denote the effective Hessian used in the untruncated Newton step
 $p^n=-(\widetilde H^n)^{-1}g^n$, and set $\delta_n:=\|\widetilde H^n-H^n\|$.
Using \Cref{lem:taylor}, we obtain, whenever \(\delta_n<\mu\),
\[
    \|e^{n+1}\|=\|e^n+p^n\|=\|(\widetilde H^n)^{-1}[(\widetilde{H}^n-H^n)e^n-r^n]\|\leq
    \frac{\delta_n}{\mu-\delta_n}\|e^n\|
    +
    \frac{M}{2(\mu-\delta_n)}\|e^n\|^2.
\]
Consequently, if $\sup_n\delta_n<\frac{\mu}{2}$, the iteration is locally linearly convergent for sufficiently small $\Vert e^0\Vert$. Superlinear convergence occurs when \(\delta_n\to0\), whereas quadratic convergence is retained if \(\delta_n=O(\|e^n\|)\). These estimates do not establish convergence
of TRHiSD-CP/TCG or variants using approximate spectral subspaces.

\end{remark}

\section{Numerical experiments}
\label{sec:4}
In this section, we measure the convergence rates of TRHiSD with the Cauchy point and Newton step in numerical examples to substantiate the theoretical results. Then, we apply TRHiSD to the LdG nematic liquid crystal model. All experiments were performed in MATLAB R2025b on a desktop computer with
an Intel Core i7-13700K CPU and 16 GB of RAM using a single computational thread.

The reported wall-clock times are averages over three runs. The eigenvectors are computed by LOBPCG using the eigenvectors from the previous outer iteration as initial approximations. The eigensolver is terminated when \(\|H^nV^n-V^n\Lambda^n\|\leq10^{-3}\), where $\Lambda^n =\diag(\lambda_1^n,...,\lambda_k^n)$. In numerical experiments, the outer iteration is terminated when the corresponding gradient residual falls below a prescribed tolerance or when the maximum iteration budget is reached. The specific stopping criteria are given in the following subsections.

\subsection{Modified Rosenbrock function}
We first consider a modified Rosenbrock function used in \cite{HiSD_Convergence}
\[
    R_m(x)=R(x)
    +
    \sum_{i=1}^d s_i\bigl[\arctan(x_i-x_i^*)\bigr]^2,
\]
where \(R(x)\) is the \(d\)-dimensional Rosenbrock function \cite{Rosenbrock1960}
\[
    R(x)
    =
    \sum_{i=1}^{d-1}
    \left[
        100(x_{i+1}-x_i^2)^2+(1-x_i)^2
    \right].
\]
Here \(x^*=(1,\ldots,1)^\top\). The additional arctangent terms are used to modify the local curvature near \(x^*\), so that \(x^*\) becomes an index-\(k\) saddle point. In the following experiments, the parameters \(s_i\) are chosen so that the Hessian of \(R_m\) at \(x^*\) has exactly \(k\) negative eigenvalues. In all Rosenbrock experiments, the iteration is terminated when the gradient residual satisfies $\|\nabla R_m(x^n)\|<10^{-5}$, or when the maximum number of iterations, \(N_{\max}=10^{5}\), is reached.

\begin{figure}[htb]
    \centering
    \includegraphics[height=5.3cm]{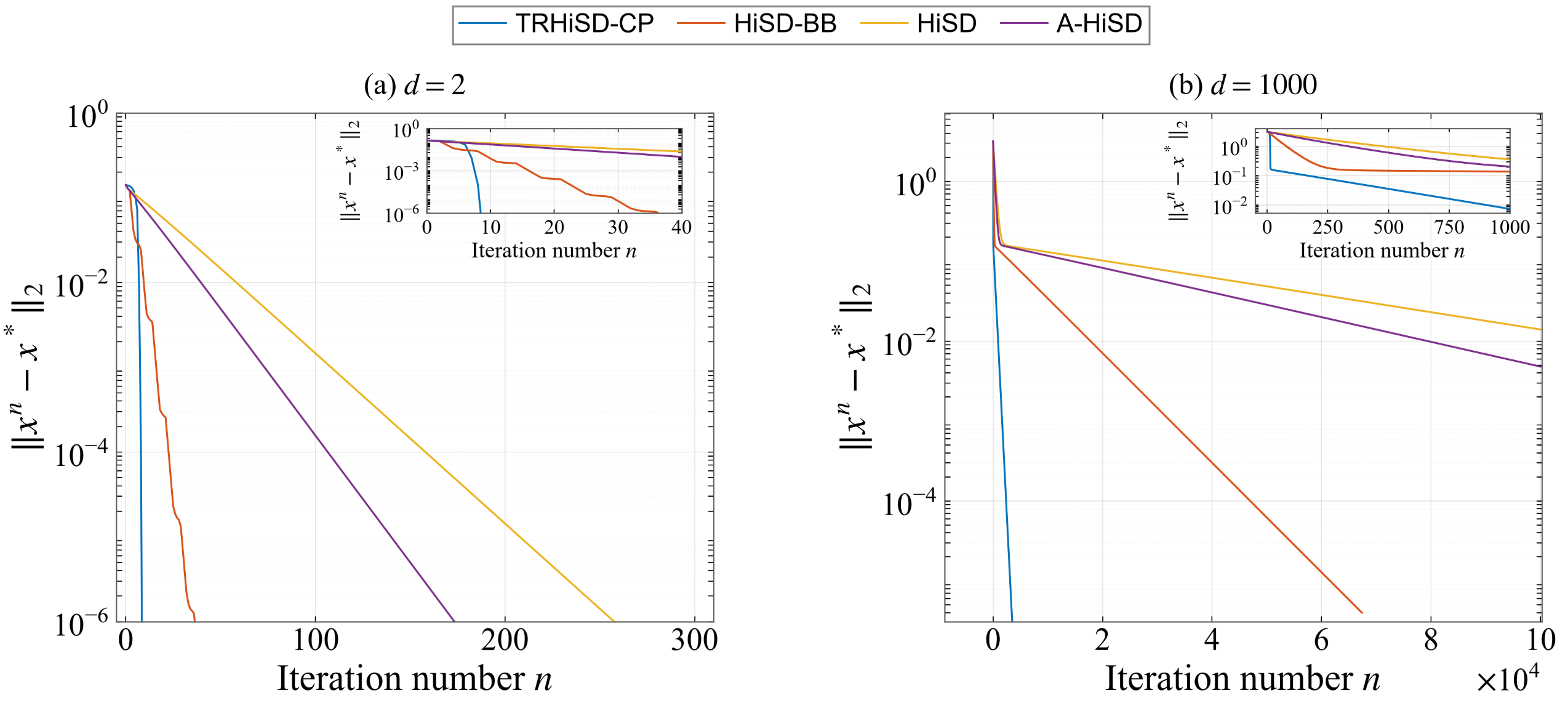}
    \caption{
    Error histories $\|x^n-x^*\|$ of TRHiSD-CP, HiSD-BB, HiSD, and A-HiSD for the modified Rosenbrock problem: (a) d=2, with the inset showing the first 40 iterations; (b) $d=1000$, with the inset showing the first 1000 iterations. The vertical axes are plotted on a logarithmic scale.
    \label{fig:cp_MR}
    }
\end{figure}

\paragraph{Efficiency comparison}
For a low-dimensional case, we take $d=2$, $s_1=-100$,  $s_2=1$.
With this choice, \(x^*=(1,1)^\top\) is an index-\(1\) saddle point. For the TRHiSD method, we set initial trust-region radii $\Delta_1^0=\Delta_2^0=10^{-3}$. 
For HiSD-BB, we choose the initial step $\beta_0=10^{-3}$, the lower bound $\beta_{\min}=\beta_0/5$ and the upper bounds $\beta_{\max}=10\beta_0$. The HiSD and A-HiSD methods use the fixed time step $\beta=10^{-3}$. \Cref{fig:cp_MR} (a) shows that TRHiSD-CP reaches the saddle point in \(9\) iterations. This corresponds to an approximately $32$-fold reduction in the number of iterations compared with the fixed-step HiSD method. HiSD-BB also converges faster than fixed-step HiSD, but TRHiSD-CP requires the fewest outer iterations in this test. 

For a higher-dimensional case, we choose $d=1000,\ s_i=-50000\ (1\leq i \leq 5),\ s_j=1\ (6\leq j)$. Then $x^*=(1,\ldots,1)^\top$ is an index-\(5\) saddle point. For the TRHiSD method, we set $\Delta_1^0=\Delta_2^0=10^{-3}$. The initial time step is chosen as $\beta_0=10^{-5}$ for HiSD-BB, while the HiSD and A-HiSD methods use the fixed time step $\beta=10^{-5}$. \Cref{fig:cp_MR} (b) shows that TRHiSD-CP also achieves faster error reduction than the fixed-step methods in this \(d=1000\) test.

\begin{figure}[htbp]
    \centering
    \includegraphics[height=4.7cm]{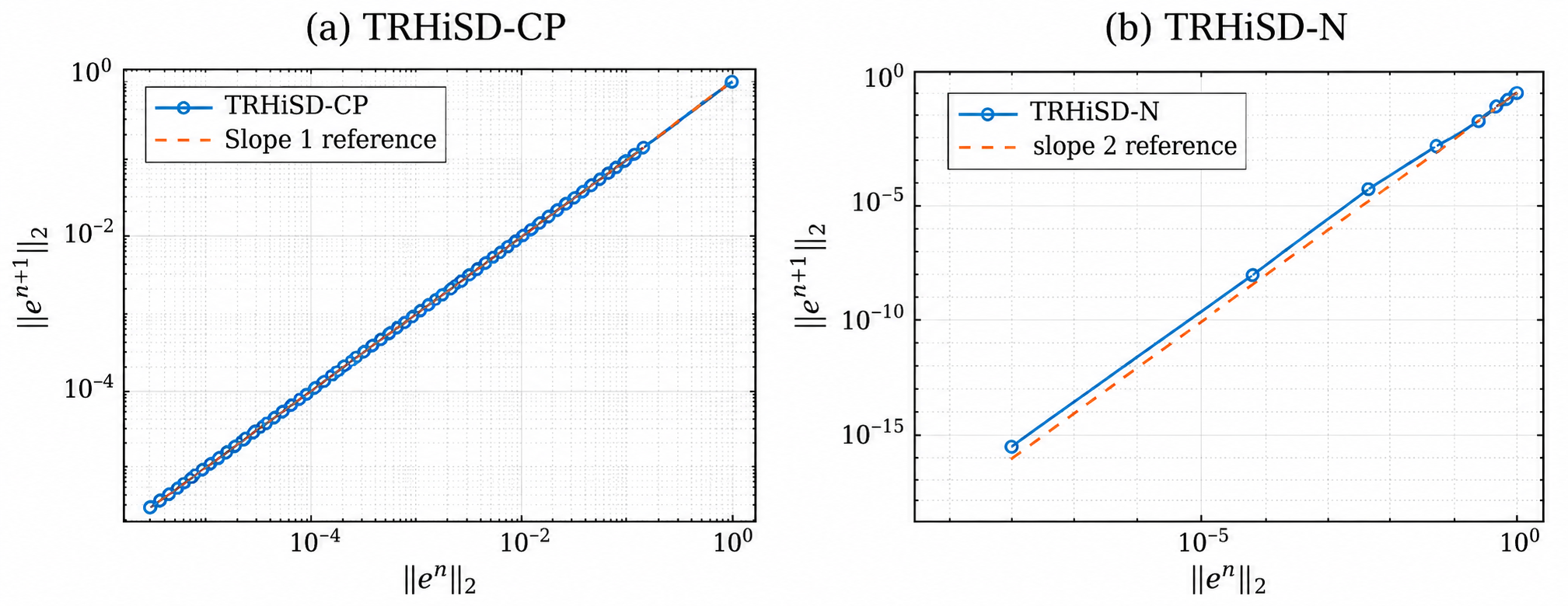}
    \caption{
    The pairs \((\|e^n\|_2,\|e^{n+1}\|_2)\) are recorded as follows: (a) every 50 iterations for TRHiSD with the Cauchy-point step; and (b) every iteration for TRHiSD with the Newton step.
    }
    \label{fig:ord_TRCP}
\end{figure}

\paragraph{Verification of local convergence rates}
We next verify the local convergence rates of the Cauchy-point and Newton-step variants of TRHiSD. We take $d=100,\ s_i=-5000\ (1\leq i \leq 3),\ s_j = 1\ (4\leq j )$. Then $x^*=(1,\ldots,1)^\top$ is an index-\(3\) saddle point. The slope in \Cref{fig:ord_TRCP} (a) is consistent with linear convergence. Over all iterations, the average contraction ratio $\frac{1}{N}\sum_{n=0}^{N-1}\frac{\|e^{n+1}\|_2}{\|e^n\|_2}$ is approximately \(0.9973\), where \(N\) is the total number of iterations, indicating slow overall linear convergence. \Cref{fig:ord_TRCP} (b) is consistent with the local quadratic convergence of the Newton-step method predicted by \Cref{thm:newton}. In the local regime, the data are consistent with $\|e^{n+1}\|_2\approx C\|e^n\|_2^2$.

\subsection{Nematic liquid crystals confined in a square}

In this experiment, we consider rodlike nematic liquid crystals confined in a two-dimensional square domain. According to the LdG theory \cite{LdGmodel}, the nematic liquid crystals can be described by a \(Q\)-tensor, which is a \(2\times 2\) symmetric traceless tensor field on the reduced square domain $ \Omega=[-1,1]^2$. 
We adopt the nondimensionalized LdG free energy \cite{PhysRevLett}
\begin{equation}
\label{eq:LdG-energy}
    E[Q]
    =
    \int_{\Omega}
    \left(
        \frac12 |\nabla Q|^2
        +
        \alpha
        \left(
            \frac{a}{4}|Q|^2
            +
            \frac18 |Q|^4
        \right)
    \right)
    \,d\mathbf r ,
\end{equation}
where \(\alpha>0\) is a parameter that describes the physical domain size. The parameter \(a\) represents the reduced temperature. Following \cite{PhysRevLett}, we set $a=-1.672$, corresponding to the nematic regime.

On the boundary \(\partial\Omega\), the nematic liquid crystals are confined to be aligned with the square boundary due to geometric confinement. Hence we impose tangent Dirichlet boundary conditions from Eqs.~(9) and (10) of the Supplemental Material in \cite{PhysRevLett}:
\begin{equation}
\label{eq:LdG-bc}
Q(\pm1,y)=-Q(x,\pm1)
=\frac{S_0}{2}
\begin{pmatrix}
1 & 0\\
0 & -1
\end{pmatrix},
\qquad x,y\in(-1,1),
\end{equation}
where \(S_0=\sqrt{2|a|}\) denotes the equilibrium scalar order parameter.

For the LdG free energy model in \eqref{eq:LdG-energy}, the solution landscape becomes increasingly complicated as the physical-size parameter \(\alpha\) increases (see Fig. 1 in the Supplemental Material of \cite{PhysRevLett}). Bifurcations may occur when the Hessian of \(E[Q]\) at a stationary solution develops a zero eigenvalue, allowing new branches of stationary solutions to emerge. 
Near bifurcation points, small-magnitude Hessian eigenvalues may occur, which can lead to a large spectral condition number and slow convergence of the Cauchy-point iteration.
More precisely, \Cref{thm:cauchy} shows that the local contraction factor of TRHiSD-CP can be arbitrarily chosen close to
\(
    q
    =
    \frac{L-\mu}{L+\mu}.
\)
In discretized LdG models, the elastic term and spatial discretization may generate a wide Hessian spectrum. Thus, when the absolute Hessian is ill-conditioned,
\(\kappa_B\gg1\), one has \(q\approx1\), and the asymptotic convergence of the Cauchy-point iteration can be very slow. To assess the efficiency of the proposed method, we employ the TRHiSD-CP/TCG variant that uses a Cauchy point on \(V^n\) and a TCG step on \(W^n\).

\begin{table}[H]
\centering
\caption{The six smallest eigenvalues of the Hessian at
selected stationary states for $\alpha=50$.
The bold entries indicate the largest negative and
smallest positive eigenvalues.}
\label{tab:eigenvalues}
\begin{tabular}{cccccccc}
\hline
State & Index & \multicolumn{6}{c}{Smallest eigenvalues} \\
\hline
$T$    & 4 & $-17.930$ & $-15.283$ & $-9.077$
       & $\mathbf{-0.415}$ & $\mathbf{2.771}$ & $5.436$ \\
$T+$   & 3 & $-17.087$ & $-14.414$
       & $\mathbf{-8.124}$ & $\mathbf{0.801}$
       & $2.817$ & $5.369$ \\
$N+$   & 4 & $-19.738$ & $-16.811$ & $-0.699$
       & $\mathbf{-0.005}$ & $\mathbf{4.373}$ & $4.896$ \\
$M+$   & 3 & $-19.732$ & $-16.815$
       & $\mathbf{-0.681}$ & $\mathbf{0.005}$
       & $4.429$ & $4.835$ \\
$BD$   & 4 & $-12.206$ & $-11.526$ & $-3.140$
       & $\mathbf{-2.420}$ & $\mathbf{8.799}$ & $8.806$ \\
$BD+$  & 3 & $-11.594$ & $-6.192$
       & $\mathbf{-2.386}$ & $\mathbf{3.167}$
       & $8.868$ & $10.441$ \\
$J$    & 2 & $-9.341$
       & $\mathbf{-1.256}$ & $\mathbf{5.067}$
       & $9.117$ & $11.576$ & $11.794$ \\
$J+$   & 1 & $\mathbf{-8.226}$ & $\mathbf{1.785}$
       & $5.367$ & $9.971$ & $11.752$ & $13.250$ \\
\hline
\end{tabular}
\end{table}

The solution landscape of the LdG model was constructed in
\cite{PhysRevLett}. We select the stationary states
\(T,\ T+,\ N+,\ M+,\ BD,\ BD+,\ J,\ J+\) and report their smallest Hessian eigenvalues in \Cref{tab:eigenvalues}. The searches $J\leftrightarrow J+,\ T\leftrightarrow T+,\ BD\leftrightarrow BD+$ are then used to compare the numerical methods.

We use a uniform \(100\times100\) grid and second-order finite differences for spatial derivatives. For the TRHiSD-CP/TCG method, we set the maximum number of TCG iterations to
\(i_{\max}\in\{20,50,100\}\) in Table 2 and the TCG tolerance to
\(\epsilon_{\mathrm{tcg}}=10^{-3}\). The upper bounds for both trust-region radii are set to \(\Delta_1^u=\Delta_2^u=10^{-3}\) for numerical stability. These upper bounds are used only as numerical safeguards and do not affect local convergence analysis. We perform numerical experiments that involve upward and downward searches between \(J\leftrightarrow J+\) , \( T+ \leftrightarrow T \) and \( BD \leftrightarrow BD+ \) at \( \alpha = 50 \). For HiSD-BB, the step size is restricted to $[2\times10^{-4},10^{-2}]$ as a numerical safeguard. For the LdG model, let \(g^n\) denote the discrete gradient vector at the \(n\)-th iterate. The outer iteration is terminated when the discrete $L^2$ norm satisfies $\sqrt{2h^2(g^n)^\top g^n}<10^{-5}$, where $h=1/50$ is the spatial mesh size, or when the maximum number of iterations, $N_{\max}=10^7$, is reached.

\Cref{tab:methods_comparison} shows that replacing the
Cauchy point in the minimization problem with a TCG step substantially reduces the number of outer iterations in all tested searches and the geometric mean time. The HiSD-BB and TRHiSD-CP require \(10^4 \sim 10^5\) iterations in several cases, while the TRHiSD-CP/TCG reduces the number of iterations to only a few hundred. These results suggest that {the Cauchy-point update in the minimization subspace leads to slow convergence in the tested cases}. By replacing this step with the TCG method, the method captures more curvature information in the minimization subspace and therefore substantially accelerates convergence. 

\begin{table}[H]
\centering
\caption{Performance comparison of HiSD-BB, TRHiSD-CP, and TRHiSD-CP/TCG for the searches between \(J\leftrightarrow J+\), \(T\leftrightarrow T+\), and \(BD\leftrightarrow BD+\) in the LdG model. The entries report the number of outer iterations. The numbers \(20,50,100\) denote the maximum number of TCG iterations used in TRHiSD-CP/TCG. The final row reports the geometric mean wall-clock time over the six searches.}
\vspace{5pt}
\label{tab:methods_comparison}
\begin{tabular}{lcccccc}
\hline
Saddle search & HiSD-BB & TRHiSD-CP & \multicolumn{3}{c}{TRHiSD-CP/TCG} \\
\cline{4-6}
& & & 20 & 50 & 100 \\
\hline
$J\to J+$ & 22888 & 81770 & 568 & 317 & 304 \\
$J+ \to J$ & 33881 & 15369 & 285 & 263 & 261 \\
$T\to T+$ & 52092 & 182396 & 518 & 231 & 200 \\
$T+ \to T$ & 79880 & 19174 & 289 & 295 & 271 \\
$BD \to BD+$ & 12280 & 39560 & 411 & 352 & 352 \\
$BD+ \to BD$ & 13950 & 9980 & 373 & 362 & 361 \\
\hline
Geometric mean time (s)
& 16.127 & 62.937 & 4.583 & 5.078 & 6.020\\
\hline
\end{tabular}
\end{table}

 We further examine the sensitivity of TRHiSD-CP/TCG to the common maximum trust-region radius and the maximum number of TCG iterations using the search $(N+\to M+)$ shown in Table 3. We take $i_{\max}\in\{25,50,75,100\}$. In this test, the upper bounds of the two trust-region radii are set equal, i.e., $\Delta_1^u=\Delta_2^u=\Delta^u$.

\begin{table}[h]
    \centering
    \small
    \caption{
        Sensitivity of TRHiSD-CP/TCG with respect to the maximum trust-region radius \(\Delta^u\) and the maximum number of TCG iterations \(i_{\max}\) for the search \(N+\rightarrow M+\). Each entry reports the number of outer iterations, with the wall-clock time in seconds given in parentheses. The symbol ``--'' indicates that the method failed to converge within the prescribed iteration budget.
    }
    \label{tab:parameter-sensitivity}
    \begin{tabular}{ccccc}
        \toprule
        \multirow{2}{*}{\(\Delta^u\)}
        & \multicolumn{4}{c}{ \(i_{\max}\) } \\
        \cmidrule(lr){2-5}
        & \(25\) & \(50\) & \(75\) & \(100\) \\
        \midrule
        \(10^{-4}\)
        & \makecell{30446 (578s)}
        & \makecell{7414 (323s)}
        & \makecell{4005 (227s)}
        & \makecell{3075 (83s)}
        \\
        \(10^{-3}\)
        & \makecell{30341 (566s)}
        & \makecell{7295 (293s)}
        & \makecell{2963 (183s)}
        & \makecell{1156 (93s)}
        \\
        \(10^{-2}\)
        & \makecell{30341 (558s)}
        & \makecell{7083 (162s)}
        & \makecell{2978 (81s)}
        & \makecell{1939 (139s)}
        \\
        \(10^{-1}\)
        & \makecell{30341 (685s)}
        & \makecell{7290 (314s)}
        & \makecell{3004 (82s)}
        & --- \\
        \bottomrule
    \end{tabular}
\end{table}

For a fixed \(\Delta^u\), increasing \(i_{\max}\) generally reduces the number of outer iterations. However, the wall-clock time does not always decrease monotonically with \(i_{\max}\). Although \(i_{\max}=100\) gives the smallest iteration count for \(\Delta^u=10^{-3}\), the smallest wall-clock time is obtained at \(\Delta^u=10^{-2}\) and \(i_{\max}=75\). 
This reflects the trade-off between reducing the number of outer iterations and increasing the potential cost of each inner TCG solver. The choice of \(\Delta^u\) is also important. If the trust-region radius is too small, for instance \(\Delta^u=10^{-4}\), the step size is overly restricted and the method may require more iterations. The parameter pair $\Delta^u=10^{-1},$ $i_{\max}=100$ did not converge within the prescribed iteration budget.

\section{Conclusions}

\label{sec:5}
We propose a TRHiSD framework for computing high-index saddle points. Under suitable assumptions and with exact spectral subspaces, we prove that the Cauchy-point variant converges locally linearly and the Newton-step variant converges locally quadratically. Numerical experiments on modified Rosenbrock functions are consistent with these convergence rates. For the LdG model, 
TRHiSD-CP/TCG, which replaces the Cauchy point in the minimization subspace with a TCG step, substantially reduces both the outer-iteration counts and the measured wall-clock time compared with HiSD-BB in the tested cases.

There are several directions for future research. The convergence analysis can be extended to accommodate inexact eigenspaces and approximate solutions of the block trust-region subproblems, while globalization strategies and global convergence can also be investigated under suitable conditions. Further improvements may involve preconditioned and quasi-Newton methods, dogleg strategies, Krylov-subspace methods, and adaptive trust-region radius updates informed by spectral information.

\bibliographystyle{siamplain}
\bibliography{references}
\end{document}